\documentclass[11pt,a4paper]{article}
\usepackage{amsfonts,amssymb}
\usepackage{calrsfs}
\usepackage[english]{babel}

\title{COMPLETENESS OF THE SPACE OF SEPARABLE MEASURES IN THE KANTOROVICH-RUBINSHTE\u{I}N METRIC}
\author{A. S. Kravchenko}
\date{} 

\begin{document}

\noindent {\bf COMPLETENESS OF THE SPACE OF SEPARABLE MEASURES IN THE KANTOROVICH-RUBINSHTE\u{I}N METRIC}
\begin{center}\bf A. S. Kravchenko\end{center}

\def\diam{{\rm diam\,}}
\def\clset{{\rm cl\,}  }
\def\id   {{\rm id\,}  }
\def\lip  {{\rm lip\,} }
\def\Lip  {{\rm Lip\,} }
\def\spt  {{\rm spt\,} }
\renewcommand{\le}{\leqslant}
\renewcommand{\ge}{\geqslant}

\renewcommand{\abstractname}{\vspace{-\baselineskip}} 
\begin{abstract}
\noindent{\bf Abstract:}
We consider the space $M(X)$ of separable measures
on the Borel $\sigma$-algebra ${\cal B}(X)$ of a metric space $X$.
The space $M(X)$ is furnished with the Kantorovich-Rubinshte\u{i}n metric known also as the
``Hutchinson distance''.
We prove that $M(X)$ is complete if and only if $X$ is complete.
We consider applications of this theorem in the theory of self-similar fractals.
\vspace{1em}

\end{abstract}

{\bf 1. Introduction.} 
In \cite{Ht}, defining the notion of an invariant measure
for a finite system of contraction similarities, Hutchinson considered a metric
space $(X,\rho)$ and the space $M_{\rm loc}(X)$ of
measures $\nu$ on $X$ with bounded support normed by the condition $\nu(X)=1$
and used the following metric on $M_{\rm loc}(X)$ \cite[4.3.(1)]{Ht}:
$$H(\nu,\mu)=\sup\left|\int_X f\,d\nu-\int_X f\,d\mu\right|,\eqno{(1)}$$
where the supremum is taken over all functions $f$ in the space
$$\lip_1(X)=\{f:X\to \mathbb{R}:\, |f(x)-f(y)|\le\rho(x,y),~\mbox{for all}~x,y\in X\}.$$

This metric, called the ``Hutchinson distance'' in \cite{Ak},
was introduced in the 1950s in the articles by
L. V. Kantorovich and G. Sh. Rubinshte\u{i}n (see \cite[3, Chapter 4, \S 4]{KnAk}).

The proof of the theorem on existence of an invariant measure
on a self-similar set in \cite[4.4(1)]{Ht} relies
upon the Banach Fixed Point Theorem;
however, the proof of completeness of $M_{\rm loc}(X)$ is absent.
In the case of a compact set $X$ the space $M_{\rm loc}(X)$
is compact (see \cite[Chapter 8, \S 4]{KnAk}) and consequently complete.
However, in general we have the following

\vspace{1ex}
 {\bf Assertion 1.1.} \sl
 If $X$ is unbounded then $M_{\rm loc}(X)$ is not complete. \rm
\vspace{1ex}

Indeed, choosing a sequence $x_k\in X$
of points such that $\rho(x_0,x_k)\le k$ and $\rho(x_0,x_k)\to\infty$ as $k\to\infty$
and using the Dirac measure $\delta_x$ (see Section 2),
we can define the sequence
$\nu_n=2^{-n}\delta_{x_0}+\sum_{k=1}^{n} 2^{-k}\delta_{x_k}$
which is Cauchy in $M_{\rm loc}(X)$ but has no limit in this space.

The gap in the proof of Hutchinson's theorem was observed in \cite{Ak}
wherein a new proof of this theorem was given in the case $X=\mathbb{R}^n$
which uses the space $M(X)$ furnished with the metric (1) of measures $\nu$
satisfying the conditions $\nu(X)=1$ and
$\int_X\rho(x_0,x)\,d\nu<\infty$ for some point $x_0\in X$.
Also, completeness of $M(X)$ was established
in the particular case $X=\mathbb{R}^n$.
In Section 4 we prove the main result (Theorem 4.2)
on equivalence of completeness of $M(X)$
and that of $X$ in the general case.
An application of this theorem to self-similar fractals (Section 5)
settles completely the question of correctness of Hutchinson's theorem
in the general case and also extends this theorem to the
case of countable systems of contractions.
This makes it possible to consider an attractor of a (countable)
system of contraction similarities
in a Banach space as the support of an invariant measure
without the a priori requirement of compactness.

This is an updated version of the original article [13]
with a fixed flow in the proof of Corollary 4.4.

\vspace{2em}
{\bf 2. Basic notions.} 
Let $(X,\rho)$ be a metric space and let $C_b(X)$ be the space of
all bounded continuous functions $f:X\to\mathbb{R}$ with the norm
$\Vert f\Vert_{\infty}=\sup_{x\in X}|f(x)|$.
The number
$$\Lip f=\sup\{|\frac{f(x)-f(y)}{\rho(x,y)}|: x\ne y\; x,y\in X\}.$$
is called the \emph{Lipschitz constant} of $f$;
if it is finite then $f$ is called a \emph{Lipschitz function}.
Denote the spaces of all real Lipschitz functions,
bounded Lipschitz functions, and functions with the Lipschitz
constant at most $\alpha$ on $X$ by
$\lip(X)$, $\lip^{\circ}(X)$, and $\lip_{\alpha}(X)$.
Introduce the following notations:
$\diam(A)=\sup\{\rho(x,y): x,y\in A\} \in[0,+\infty]$
is the diameter of a set $A\subset X$ and
$\rho(x,A)=\inf_{y\in A}\rho(x,y)$
is the distance from a point $x\in X$ to a set $A$.
The \emph{support of a real function} $f:X\to [0,+\infty)$ is the set
$\spt f=\{x\in X: f(x)>0\}$.

We write $\alpha_n\downarrow\alpha_0$ ($\alpha_n\uparrow\alpha_0$)
for a monotonically decreasing (increasing) numeric
sequence $\{\alpha_n\}$ convergent to $\alpha_0$.
Similarly, we write ${f_n\downarrow f_0}$ (${f_n\uparrow f_0}$)
on $X$ if the sequence $\{f_n\}$ of real functions on $X$
decreases (increases) and converges pointwise to $f_0$ everywhere on $X$.

By a \emph{measure} $\nu$ on $X$ we mean a nonnegative
countably additive real set function given on the $\sigma$-algebra
${\cal B}(X)$ of all Borel subsets of $X$
which satisfies the equality $\nu(\varnothing)=0$.
The measure
$\delta_x(A)=\{1 ~\mbox{for}~ x\in A;~ 0 ~\mbox{for}~ x\not\in A\}$
is called the \emph{Dirac measure} at $x\in X$ (see [4, 10.9.4(1)]).
A measure $\mu$ on $X$
satisfying the condition $\mu(X)<+\infty$ is said to be \emph{finite}.
The \emph{support} of a measure $\mu$ is
$\spt\mu=X\setminus\cup\{A\subset X: A ~\mbox{is open and}~\mu(A)=0 \}.$
A set $A\subset X$ is called separable if it is contained in the closure
of an at most countable subset.
In line with \cite[Chapter 1,\S 1.]{Bl},
we say that a measure $\mu$ is \emph{separable} if there is
a separable Borel set $A\subset X$ such that $\mu(X\setminus A)=0$.
It is well known \cite[2.2.16]{Fd} that the support of
a finite measure is always separable and therefore for finite measures the condition $\mu(X\setminus\spt\mu)=0$
is equivalent to the separability of $\mu$.

A family $\Pi$ of finite measures is called \emph{dense} (see [5, Chapter 1, \S 1]) if, for every $\varepsilon>0$, there
is a compact set $K$ such that $\nu (X\setminus K)<\varepsilon$ for all measures $\nu\in\Pi$;
respectively, a finite measure is called \emph{dense}
if so is the family constituted by this single measure. A finite dense measure is determined
uniquely by the values of the integral at $f\in\lip^{\circ}(X)$
(see \cite[5, Chapter 1, \S 1, Theorem 1.3]{Bl}).
Observe that, in a complete metric space $X$, the density
and separability properties of a measure are equivalent
(see \cite[Supplement III]{Bl}).
We only consider separable measures, omitting the so-called
measure problem (see \cite[Supplement III]{Bl}) or \cite[2.1.6]{Fd}).
Each finite measure on $X$ generates a linear functional
$\mu(f)=\int_X f\,d\mu$ on the space $E(X,\mu)$
of all $\mu$-summable real functions.
In line with \cite[p. 159]{Ak}, given a metric space $X$
we consider the space $M(X)$ of all separable measures $\mu$ such that $\mu(X)=1$ and $\mu(f)<+\infty$ for every
function $f\in\lip(X)$. Note that in general $M_{loc}(X)\subset {\cal M}(X)$,
where ${\cal M}(X)$ is the space of measures
without the separability requirement as introduced in \cite[p. 160]{Ak}.
Moreover, ${M(X)\subset {\cal M}(X)}$.
With each point $a\in X$ we associate the function
${\phi_a(x)=\rho(a,x)}$ in $\lip(X)$.

\vspace{1em}

{\bf Assertion 2.1.} \sl 
 A separable measure $\mu$ on $(X,\rho)$, satisfying the condition $\mu(X)=1$, belongs
to $M(X)$ if and only if there is a point $a\in X$ such that $\mu(\phi_{a})<+\infty$. \rm

{\tt Proof.}
If $\mu\in M(X)$ then $\phi_a\in\lip(X)$ for every point $a\in X$ and consequently $\mu(\phi_{a})<+\infty$.
If there is a point $a\in X$ such that $\mu(\phi_{a})<+\infty$
then the estimate
$$f(x)\le f(a)+(\Lip f)\cdot\rho(a,x)=f(a)+(\Lip f)\cdot \phi_a(x)$$
holds for every function $f\in\lip(X)$,
which yields the sought estimate for the integral:
$\mu(f)\le f(a)+(\Lip f)\cdot \mu(\phi_a)< +\infty$.
\hbox{$\Box$}

\vspace{2em}
{\bf 3. Metrization of the space of measures.} 
The distance $H(\mu,\nu)$ between the measures $\mu$ and $\nu$ is defined by (1),
where the supremum is taken over all functions $f\in\lip_1(X)$,
and is the restriction to $M(X)$ of the metric $H(\mu,\nu)$ considered in [1]
on the formally broader space ${\cal M}(X)$.
Note that the completeness property of $X$ was not used in
checking the axioms of a metric for $H(\mu,\nu)$ in \cite[Theorem 1, p. 161]{Ak}.
We say that a sequence $\mu_k$ of finite measures
{\it converges weakly to a finite measure} $\mu$
if $\mu_k(f)\to\mu(f)$ as $k\to\infty$ for every function $f\in C_b(X)$;
in this case we write $\mu_k\Rightarrow\mu$. A family of sets of the form
$\{\mu\in M(X): |\mu(f_j)-\nu(f_j)|<\varepsilon; j=1,\ldots,k\}$
for arbitrary $\varepsilon>0$ and arbitrary finite collections
$f_1,\ldots,f_k$
of functions in $C_b(X)$ determines a system of basic neighborhoods
for each point $\nu\in M(X)$
which generates the topology of {\it weak convergence} which we denote by ${\cal W}$ .

A sequence $\{\nu_n\}_{n=1}^{\infty}$
is called {\it weakly Cauchy} if $\nu_n(f)$ is Cauchy for every function
$f\in C_b(X)$. A set $\Pi$ of measures is {\it weakly complete} if
every weakly Cauchy sequence $\{\nu_n\}$ in this set converges weakly to
some measure in $\Pi$.


\filbreak 
{\bf Theorem 3.1.} \sl
For every metric space $X$, the topology ${\cal T}$
on $M(X)$ generated by $H(\mu,\nu)$
coincides with the topology ${\cal W}$ of weak convergence if and only
if $\diam X<\infty$. Moreover, if $\diam(X)=\infty$
then ${\cal T}$ is strictly finer than ${\cal W}$. \rm
\vspace{1em}

Before proving this theorem, we give several auxiliary assertions:

\vspace{1em} {\bf Lemma 3.2.} \sl 
For every function $f\in C_b(X)$,
there is a sequence $\{\varphi_n\}$ of functions in
$\lip^{\circ}(X)$ such that $\varphi_n\downarrow f$ on $X$. \rm
\vspace{1em}

{\tt Proof.} 
Let $f\in C_b(X)$ and $||f||_{\infty}=m<+\infty$.
For $n=1,2,\ldots$ each of the continuous functions
$\varphi_n(x)=\sup\{f(t)-n\cdot\rho(x,t): t\in X\}$
is bounded on $X$, since
$-m\le f(x)\le\varphi_n(x)\le\sup\{f(t): t\in X\}=m.$
It follows from the inequality
$(f(t)-n\cdot\rho(x_1,t))-n\cdot\rho(x_1,x_2) \le f(t)-n\cdot\rho(x_2,t) \le \varphi_n(x_2)$
that
$\varphi_n(x_1)-\varphi_n(x_2)\le n\cdot\rho(x_1,x_2)$
for arbitrary $x_1,x_2\in X$.
Consequently, $\Lip\varphi_n\le n$ and $\varphi_n\in\lip^{\circ}(X).$
The inequality ${\varphi_{n+1}(x)\le\varphi_n(x)}$
is immediate from the definition of the functions $\varphi_n(x)$.
For every fixed $x\in X$ and an arbitrary $\varepsilon>0$, there is
$\delta=\delta(\varepsilon)>0$ such that $|f(x)-f(t)|<\varepsilon$
whenever $\rho(x,t)<\delta$.
Since the inequality
$(f(t)-f(x))-n\cdot\rho(x,t)\le 0$
holds for all $n>2m/\delta$ and $\rho(x,t)\ge\delta$, the estimate
$$0\le\varphi_n(x)-f(x)=\sup\{(f(t)-f(x))-n\cdot\rho(x,t): \rho(x,t)<\delta\}\le\varepsilon$$
is valid for all sufficiently large n. Consequently,
$f(x)=\lim_{n\to\infty}\varphi_n(x)$. \hbox{$\Box$}
\vspace{1em}

{\tt Remark.}
A similar assertion for semicontinuous functions
on bounded sets in $\mathbb{R}^n$ was proven by
Hausdorff (1919); see \cite[Theorem II.5]{Ts}.
\vspace{1em}

Let $V$ be a vector lattice. A linear functional
$F:V\to\mathbb{R}$ is {sequentially $o$-continuous}
or {sequentially order continuous}
if $F(u_n)\to 0$ as $n\to\infty$ for every monotonically decreasing sequence $u_n$ of elements
of $V$ such that $\inf u_n=0$.
An example of a sequentially $o$-continuous functional is the integral
with respect to an arbitrary finite measure $\mu$ on $X$.

\vspace{1em} {\bf Lemma 3.3.} \sl
Suppose that $V$ is a vector lattice and $V_0$ is a subset of $V$ such that,
for every element $v\in V$,
there is a decreasing sequence $u_k\in V_0$, $k=1,2,\ldots$, such that $\inf u_k=v$.
If a sequence
$F_n:V\to\mathbb{R}$, $n=0,1,\ldots$
of sequentially $o$-continuous positive linear functionals
converges pointwise to $F_0$ on $V_0$
then the sequence $\{F_n\}$ converges pointwise to $F_0$
on the whole $V$. \rm

\vspace{1em} {\tt Proof.}
Choose $v\in V$ and $\varepsilon>0$.
By condition, there is a decreasing sequence $u_k\in V_0$
for which $\inf u_k=v$.
Choose $k=k(\varepsilon)$ such that ${|F(u_k)-F(v)|<\varepsilon}$
and $N=N(\varepsilon)$ such that ${|F_n(u_k)-F(u_k)|<\varepsilon}$
for ${n>N}$. Then
$$F_n(v)\le F_n(u_k)< F(u_k)+\varepsilon< F(v)+2\varepsilon.$$
Passing to $\limsup$ as $n\to\infty$, we obtain
$\limsup F_n(v)\le F(v)+2\varepsilon$.
In view of the arbitrariness of $\varepsilon$, we have lim sup $\limsup F_n(v)\le F(v)$.
Applying similar arguments to $(-v)$, we find that
${\liminf F_n(v)\ge F(v)}$ whence ${\lim F_n(v)=F(v)}$.
\hbox{$\Box$}

\vspace{1em} {\bf Corollary 3.4.} \sl
Suppose that a sequence $\{\mu_n\}$ of finite measures
and a finite measure $\mu$ are such
that $\mu_n(f)\to\mu(f)$ for every function $f\in\lip^{\circ}(X)$
Then $\mu_n\Rightarrow\mu$. \rm
\vspace{1em} 

{\tt Proof.}
Put $V=C_b(X)$ and $V_0=\lip^{\circ}(X)\subset V$.
Then the result of the corollary is immediate from Lemmas 3.2 and 3.3.
\hbox{$\Box$} \vspace{1em}

{\bf Corollary 3.5.} \sl
The topology ${\cal T}$ on $M(X)$ generated by $H$
is not coarser than ${\cal W}$. \rm
\vspace{1em} 

{\tt Proof.}
It follows from Corollary 3.4 that every sequence in $M(X)$ convergent in the metric H
converges in the weak topology.
\hbox{$\Box$}

\vspace{1em} {\bf Lemma 3.6.} \sl
If $X$ is a bounded space then every weakly convergent sequence in $M(X)$ converges
in the metric H. \rm
\vspace{1em} 

{\tt Proof.}
Assume that $X$ is bounded.
Consider an arbitrary weakly convergent sequence
$\nu_n\Rightarrow\nu_0$ of separable measures ($n=1,2,\ldots$).
Let $\tilde{X}$ be the completion of $X$. The family of sets
${{\cal D}=\{A\subset\tilde{X}: A\cap X\in{\cal B}(X)\}}$
is a $\sigma$-algebra.
This is immediate from the fact that the family ${\cal B}(X)$ of Borel sets is a $\sigma$-algebra.
If a set V is open in $\tilde{X}$
then $V\cap X$ is open in $X$; consequently,
$V\cap X\in{\cal B}(X)$ and $V\in{\cal D}$.
Since ${\cal B}(\tilde{X})$ is the minimal $\sigma$-algebra
containing all open sets in $\tilde{X}$ , we have
${\cal B}(\tilde{X})\subset{\cal D}$ or, in other words,
$\{ A\cap X: A\in{\cal B}(\tilde{X})\}\subseteq {\cal B}(X).$
Thus, each separable measure $\nu_n:{\cal B}(X)\to[0,+\infty)$
generates a separable measure $\tilde{\nu}_n:{\cal B}(\tilde{X})\to[0,+\infty)$
defined by the formula
$\tilde{\nu}_n(A)=\nu_n(A\cap X)$ 
(all necessary properties of $\tilde{\nu}$ are immediate from the similar
properties of the measure $\nu$).

For every function $f\in C_b(\tilde{X})$ we have
$$\tilde{\nu}_n(f)=\nu_n(f|_X)\to\nu_0(f|_X)=\tilde{\nu}_0(f),$$
i.e., $\tilde{\nu}_n\Rightarrow\tilde{\nu}_0$.
Consequently, the family $\Pi=\{\tilde{\nu}_n\}_{n=0}^{\infty}$
is compact in the topology of weak convergence.
Fix $\varepsilon>0$.
By the converse Prokhorov theorem \cite[Chapter 1, \S 6, Theorem 6.2]{Bl},
the family $\Pi$ is dense; i.e.,
there is a compact set $K\subset \tilde{X}$ such that
$\tilde{\nu}_n(\tilde{X}\setminus K)<\varepsilon$ for $n=0,1,\ldots$.

Since the set $X$ is everywhere dense in $\tilde{X}$,
the collection of balls $B(x,\varepsilon)$
of diameter $\varepsilon$ centered at $x\in X$
covers $\tilde{X}$ and in particular the compact set $K$.
Choose a finite subcovering of $K$:
$\cup_{i=1}^{m} B(x_i,\varepsilon)\supset K$, and denote
$A=\{x_i\}_{i=1}^{m}$.
Since the set $A$ is finite, the set of functions
$$\mathcal{F}=\{\varphi\in\lip_1(A): \varphi(x_1)=0 \}$$
is homeomorphic to a closed bounded subset of $\mathbb{R}^m$
and consequently $\mathcal{F}$ is compact in the uniform norm
$\|\varphi\|_{\infty}=\max|\varphi|.$
Choose a finite $\varepsilon$-net
$\varphi_k\in\mathcal{F}$ for $\mathcal{F}$.
Consider the collection of functions $\psi_k:X\to\mathbb{R}$:
$$\psi_k(x)=\min\{\varphi_k(x_i)+\rho(x,x_i): i=1,\ldots,m\}.$$
Each of these functions is an extension of $\varphi_k$
to the whole space $X$,
${\psi_k|_A\equiv\varphi_k}$;
moreover, $\psi_k\in\lip_1(X)$ by construction.

Choose a number $n_0$ so large that
$|\nu_n(\psi_k)-\nu_0(\psi_k)|<\varepsilon$
for all $k$ and all $n\ge n_0$.
Let $f\in\lip_1(X)$ be an arbitrary Lipschitz function.
Then the function $g(x)=f(x)-f(x_1)$ is such that $g|_A\in\mathcal{F}$
and hence there is $k$ with the property
$\max\{|g(x_i)-\varphi_k(x_i)|: i=1,\ldots,m\}<\varepsilon$.
For every point $x\in K\cap X$, there is a point $x_i$ such that
$\rho(x_i,x)\le\varepsilon$ and then
$$|g(x)-\psi_k(x)|\le|g(x)-g(x_i)|+|g(x_i)-\psi_k(x_i)|+|\psi_k(x_i)-\psi_k(x)|<3\varepsilon.$$
Thus,
$\max\{|g(x)-\psi_k(x)|: x\in K\cap X\}\le 3\varepsilon$.
Gathering the above estimates, for $n\ge n_0$ we obtain
$$|\nu_n(f)-\nu_0(f)|=|\nu_n(g)-\nu_0(g)|\le
  \biggl|\int_{X\cap K} g\,d(\nu_n-\nu_0)\biggr|+\biggl|\int_{X\setminus K} g\,d(\nu_n-\nu_0)\biggr|\le$$
$$\le|\nu_n(\psi_k)-\nu_0(\psi_k)| + 
  2\max_{x\in X\cap K}|g(x)-\psi_k(x)| +
  \max|g|\cdot|\nu_n(X\setminus K)+\nu_0(X\setminus K)|\le$$
$$\le\varepsilon + 6\varepsilon + 2\varepsilon\cdot\max|g| \le
  \varepsilon\cdot(7+2\diam X)$$
Passing to the least upper bound over the functions $f\in\lip_1(X)$, we arrive at the estimate
$$H(\nu_n,\nu_0)\le \varepsilon\cdot(7+2\diam X).$$
In view of the arbitrariness of $\varepsilon$,
this proves convergence of $\{\nu_n\}$ to $\nu_0$ in the metric $H$.
\hbox{$\Box$} 

\vspace{1em} {\bf Lemma 3.7.} \sl
 If $X$ is an unbounded space then the topology of the metric $H$ on $M(X)$ is not
equivalent to (and is strictly finer than) the weak topology. \rm
\vspace{1em} 

{\tt Proof.}
Assume that $X$ is unbounded. Choose a sequence $\{x_n\}_{n=0}^{\infty}$
n=0 of points such that $\rho(x_0,x_n)\ge n^2$.
Let
$\nu_n = \frac{1}{n}\cdot\delta_{x_n} + (1-\frac{1}{n})\cdot\delta_{x_0},$
Then $\nu_n \Rightarrow \delta_{x_0}$.
Indeed, for every bounded continuous function $f$ we have
$$|\delta_{x_0}(f)-\nu_n(f)| =
  \frac{1}{n}\cdot |f(x_0)-f(x_n)| \le
  \frac{1}{n}\cdot 2\max|f|\mathop{\longrightarrow}\limits_{n\to\infty} 0.$$
On the other hand,
$$H(\delta_{x_0},\nu_n)=\nu_n(\phi_{x_0}) = \frac{1}{n}\cdot\rho(x_0,x_n) \ge n,$$
and so convergence in the metric H is absent.
\hbox{$\Box$} \vspace{1em}

Theorem 3.1 is a straightforward consequence of 3.5--3.7.

\vspace{2em}
{\bf 4. The basic theorems on completeness of the space of measures.} 
In this section we prove the following two theorems:

\vspace{1em} {\bf Theorem 4.1.} \sl
If $X$ is a complete space then the space
of finite separable measures on $X$ is weakly complete. \rm

\vspace{1em} {\bf Theorem 4.2.} \sl
The space $M(X)$ is complete in the metric $H$
if and only if $X$ is complete. \rm
\vspace{1em} 

We start with technical lemmas. The following lemma appears in [5, Chapter 1, \S 6] as a part of the
proof of the converse assertion of the Prokhorov theorem.

\vspace{1em} {\bf Lemma 4.3.} \sl
In a complete space $X$, a family $\Pi$ of measures is dense if and only if, for arbitrary
$\varepsilon>0$ and $\delta>0$, there is
a finite collection $\{B_i\}$ of balls of radius $\varepsilon$ such that
$\nu(X\setminus \cup B_i)<\delta$ for all measures $\nu\in\Pi$. \rm
\vspace{1em} 

The proof of necessity is immediate from the definition of
a dense family and total boundedness of
a compact set.
Conversely, let $U_k=\cup_i B_{ik}$
be the union of a finite collection of balls of radius $\frac{1}{2^k}$
and let $\nu(X\setminus U_k)<\frac{\delta}{2^k}$.
Then the intersection of the sets $U_k$
is totally bounded by construction. Hence, by
Hausdorff's theorem {\cite[4.6.7.]{Kut}}, its closure $K=\clset\cap_k{U_k}$
is compact and moreover $\nu(X\setminus K)<\delta$.
\hbox{$\Box$} \vspace{1em}

The Lemma 4.3 directly implies the foollowing

\vspace{1em} {\bf Corollary 4.3.1.} \sl
In a complete space $X$, for a not dence family $\Pi$ of measures
there exist $\varepsilon>0$ and $\delta>0$ such that,
for every finite collection $\{B_i\}$ of balls of radius $\varepsilon$,
there is a measure $\nu\in \Pi$ such that
$\nu(X \setminus \cup B_i) \ge \delta$. \rm

\vspace{1em} {\bf Corollary 4.4.} \sl
If a sequence $\{\nu_k\}$ of separable measures is not dense
then there exist $\varepsilon>0$ and $\delta>0$ such that,
for every finite collection $\{B_i\}$ of balls of radius $\varepsilon$
and every $n_0\ge 0$, there is a natural number $n>n_0$
such that $\nu_n(X \setminus \cup B_i) \ge \delta$. \rm
\vspace{1em}

{\tt Proof.}
Assume that a sequence
$\Pi=\{\nu_k\}_{k=1}^\infty$
of separable measures is not dense.
There exist $\varepsilon>0$ and $\delta>0$ such that,
the assertion of the Corollary 4.3.1 holds.

Choose an arbitrary $n_0>0$ and
a finite collection $\{B_i\}$ of balls of radius $\varepsilon$.
The finite family $\Pi'=\{\nu_k: k\le n_0\}$ is dense.
By applying Lemma 4.3 to the choosen $\varepsilon$, $\delta$ and $\Pi'$
we obtain that there exists a finite collection $\{B'_i\}$
of balls of radius $\varepsilon$ such that
$\nu(X\setminus \cup B'_i)<\delta$ for all measures $\nu\in\Pi'$.
Since adding more balls to $\{B'_i\}$ will not break this inequality,
we can assume that $\{B'_i\}$ collection includes $\{B_i\}$.

From the Corollary 4.3.1 and the selection of $\varepsilon$ and $\delta$ it follows,
that there exists $\nu_n\in \Pi$ such that
$\nu_n(X \setminus \cup B'_i) \ge \delta$.
That means $\nu_n\not\in\Pi'$ and $n>n_0$.
The inclusion $\cup B_i\subset\cup B'_i$ implies the desired inequality
$\nu_n(X \setminus \cup B_i) \ge \delta$.
\hbox{$\Box$}

\vspace{1em} {\bf Lemma 4.5.} \sl
Let $X$ be a complete space.
If a sequence $\{\nu_n\}$ of finite separable measures is such
that, for every function $f\in\lip^{\circ}(X)$,
$\nu_n(f)$ is Cauchy
then $\{\nu_n\}$ is a dense family.
As a consequence, every weakly Cauchy sequence
of separable measures is a dense family. \rm
\vspace{1em} 

{\tt Proof.} 
Assume that the sequence $\{\nu_n\}$ is not dense. Show that there is
a function $f\in\lip^{\circ}(X)$
such that $\nu_n(f)$ is not Cauchy. Denote by
${B(x,r)=\{t\in X: \rho(t,x)\le r\}}$ the closed ball
of radius $r>0$ centered at $x\in X$.
Given a finite set $A$ and $r>0$, we denote
$A^r=\cup_{x\in A} B(x,r)$.
It follows from
Corollary 4.4 that there exist $\varepsilon>0$ and $\delta>0$ such that,
for every finite set $A\subset X$ and every $n_0\ge 0$ ,
there is a natural number $n>n_0$ such that
$\nu_n(X \setminus A^{\varepsilon}) \ge \delta$.

Construct a sequence of natural numbers $1\le n_1<n_2<\ldots$
and sequences $\{A_k\}$ and $\{D_k\}$ of finite
subsets of $X$ such that the following properties are satisfied:

$\begin{array}{ll}
(i) & A_{i-1}\subset A_i, ~\mbox{for all}~ i>1 \\
(ii) & D_i\subset A_j, ~\mbox{for}~ i\le j \\
(iii) & D_i^{\varepsilon/2} \cap A_j^{\varepsilon/2} = \emptyset, ~\mbox{for}~ i>j \\
(iv) & \nu_{n_k} \bigl( D_k^{\varepsilon/4} \bigr) >\delta/2 \\
(v) & \nu_{n_k} \bigl( X \setminus A_k^{\varepsilon/2} \bigr) <\delta/32
\end{array}$

We proceed by the following algorithm.

{\tt Step 1.} By the choice of $\varepsilon$ and $\delta$,
there is a number $n_1$ such that $\nu_{n_1}(X)\ge\delta$.
By separability of $\nu_{n_1}$,
there exists a finite set $A_1$ such that
$\nu_{n_1} \bigl( X\setminus A_1^{\varepsilon/4} \bigr) <\delta/32$
Property $(v)$ for $k=1$ follows from the embedding
$A_1^{\varepsilon/4}\subset A_1^{\varepsilon/2}$
and monotonicity of the measure $\nu_{n_1}$.
Put $D_1=A_1$. Then $(ii)$ is valid for $i=j=1$ and
$\nu_{n_1} \bigl( D_1^{\varepsilon/4} \bigr) =
 \nu_{n_1}(X) - \nu_{n_1} \bigl( X\setminus A_1^{\varepsilon/4} \bigr)
 > \delta - \delta/32 > \delta/2$,
which is nothing but $(iv)$ for $k=1$.

{\tt Step k for $k>1$.}
Assume there are numbers
$\{n_i\}$ and sets $A_i$ and $D_i$, $i=1,\ldots,k-1$,
with properties $(i)-(v)$.
It follows from the choice of $\varepsilon$ and $\delta$
that there is a number $n_k>n_{k-1}$ such that
$\nu_{n_k} \bigl( X\setminus A_{k-1}^{\varepsilon} \bigr) \ge \delta$.
Since $\nu_{n_k}$ is separable then there is a
finite set
$D_k\subset X\setminus A_{k-1}^{\varepsilon}$
with property $(iv)$. Since
$D_k\cap A_{k-1}^{\varepsilon}=\emptyset$, we have
$D_k^{\varepsilon/2}\cap A_{k-1}^{\varepsilon/2}=\emptyset$.
From $(i)$ we deduce $A_j\subset A_{k-1}$ for $j<k-1$, which
implies $(iii)$ in the case $i=k$.
By separability of $\nu_{n_k}$, there is a finite set $F_k$ such that
$\nu_{n_k} \bigl( X\setminus F_k^{\varepsilon/2} \bigr) < \delta/32$.
Put $A_k=F_k\cup A_{k-1}\cup D_k$, whence we obtain
$(i)$ in the case $i=k$, $(ii)$ in the case $j=k$, and $(v)$.

Take the following sequences $\{\varphi_k\}$ and $\{f_k\}$ of real functions:
$$\varphi_k(x) = \max(1-\rho(x, D_k) \cdot 2 / \varepsilon, 0), ~\mbox{where}~k=1,2,\ldots,$$
$$f_1=0,~f_{k+1}=
  \left\{
  \begin{array}{ll} 
   f_k & \mbox{if}~|\nu_{n_k}(f_k)-\nu_{n_{k+1}}(f_k)|>\frac{\delta}{8}, \\
   f_k+\varphi_{k+1} & \mbox{otherwise}.
  \end{array}
  \right.$$

It follows from $(i)$--$(iii)$ that
$D_i^{\varepsilon/2}\cap D_j^{\varepsilon/2}=\emptyset$ for $i\ne j$.
Thus, the supports of the functions
$\spt\varphi_k=D_k^{\varepsilon/2}$, $k\in \mathbb{N}$
are pairwise disjoint.
For every $k\in\mathbb{N}$ and $x\in D_k^{\varepsilon/4}$ we have
$\varphi_k(x)\ge 1/2$. Hence,
$\nu_{n_k}(\varphi_k) \ge \frac{1}{2} \nu_{n_k} \bigl( D_k^{\varepsilon/4} \bigr)$.
Applying $(iv)$, we obtain $\nu_{n_k}(\varphi_k)>\frac{\delta}{4}$
which implies that the estimate
$$|\nu_{n_k}(f_k)-\nu_{n_{k+1}}(f_{k+1})|>\delta/8$$
holds for every $k\in\mathbb{N}$.
Define $f(x)=\sup_k f_k(x)$.
Since $0\le\varphi_k(x)\le 1$, $\Lip \varphi_k \le 2/\varepsilon$,
and the supports of the functions $\varphi_k$ are disjoint, we have
$0\le f_k\le 1$ and $\Lip f_k \le 2/\varepsilon$ for all $k\in\mathbb{N}$;
moreover, $0\le f\le 1$ and $\Lip f\le 2/\varepsilon$.
Thus, $f\in\lip^{\circ}(X)$.
For every natural $k$
$$0\le (f-f_k)\le \sum_{i=k+1}^{\infty}\varphi_i \le 1.$$
Consequently,
$$\spt (f-f_k)\subset \cup_{i=k+1}^{\infty} D_i^{\varepsilon/2}.$$
By $(iii)$,
$\cup_{i=k+1}^{\infty} D_i^{\varepsilon/2}\subset X\setminus A_k^{\varepsilon/2}$
and
$$|\nu_{n_k}(f-f_k)|\le \nu_{n_k} \bigl( X\setminus A_k^{\varepsilon/2} \bigr) < \varepsilon/32.$$
Applying the triangle inequality, we obtain
$$|\nu_{n_k}(f)-\nu_{n_{k+1}}(f)| \ge |\nu_{n_k}(f_k)-\nu_{n_{k+1}}(f_{k+1})| - $$
$$ - |\nu_{n_k}(f)-\nu_{n_k}(f_k)|-|\nu_{n_{k+1}}(f)-\nu_{n_{k+1}}(f_{k+1})| >
\frac{\delta}{8}-2\frac{\delta}{32}=\frac{\delta}{16}.$$
Hence, $\nu_n(f)$ is not Cauchy.
\hbox{$\Box$}

\vspace{1em} {\bf Lemma 4.6.} \sl
Suppose that $X$ is a complete space and a sequence $\{\nu_n\}$ of finite separable measures is
such that, for every function $f\in\lip^{\circ}(X)$,
the numeric sequence $\nu_n(f)$ is Cauchy. Then there is a unique
separable measure $\mu$ which the sequence $\{\nu_n\}$ converges weakly to. \rm
\vspace{1em} 

{\tt Proof.} 
\sloppy Define a nonnegative linear functional $\mu$ by the formula
${\mu(f)=\lim\limits_{n\to\infty}\nu_n(f)}$.
for all functions $f\in\lip^{\circ}(X)$.
Consider a monotonically vanishing sequence $\varphi_k\in\lip^{\circ}(X)$.
Fix $\varepsilon>0$.
By Lemma 4.5,
the sequence $\nu_n$ constitutes a dense family;
hence, there is a compact set $K$ such that
$\nu_n(X\setminus K)<\varepsilon$.
The convergent sequence $\varphi_k$ converges
uniformly on a compact set; consequently, there is a number $k_0$
so large that
${\max\{\varphi_k(x): x\in K\}<\varepsilon}$ for $k>k_0$ and
$$\nu_n(\varphi_k)=\int_K \varphi_k\,d\nu_n +
  \int_{X\setminus K} \varphi_1\,d\nu_n \le
  \varepsilon \cdot(\nu_n(1) + \max\varphi_1),$$
$$\mu(\varphi_k)=\lim\limits_{n\to\infty}\nu_n(\varphi_k)\le 
  \varepsilon \cdot(\mu(1) + \max\varphi_1)$$
for all $k>k_0$.
Thus, $\mu(\varphi_k)\downarrow 0$
for an arbitrary sequence $\varphi_k \downarrow 0$.
By Daniel's theorem \cite[II.7.1]{Nv},
the functional $\mu$ is a measure on $X$.

To show the separability of $\mu$, it suffices to consider
a separable set $S={\rm cl}(\cup\, \spt\nu_n)$ and demonstrate
that $\mu(S)=\mu(X)$.
Using a sequence of the Lipschitz functions
$f_k(x)=\max( 1-k\cdot\rho(S,x), 0)$,
which is monotonically convergent everywhere on $S$
to the characteristic function $\chi_S$ of $S$, we obtain
$$\mu(S)=\mu(\chi_S)=\lim\limits_{k\to\infty}\mu(f_k)=
  \lim_{k\to\infty}\lim_{n\to\infty}\nu_n(f_k)=\mu(1)=\mu(X).$$
Applying Corollary 3.4, we conclude that the sequence
$\{\nu_n\}$ converges weakly to a separable measure $\mu$.
\hbox{$\Box$} \vspace{1em}

{\tt Proof of Theorem 4.1.}
 Consider an arbitrary weakly Cauchy sequence $\{\nu_n\}$ of separable measures,
i.e., such that for every function ${f\in C_b(X)}$ the numeric sequence
$\{\nu_n(f)\}$ is Cauchy. Since
$\lip^{\circ}(X)\subset C_b(X)$, by Lemma 4.6, the sequence
$\{\nu_n\}$ converges weakly to some separable measure.
Hence, the space of separable measures is weakly complete.
\hbox{$\Box$} \vspace{1em}

{\tt Proof of Theorem 4.2.} 
Assume that a sequence $\{\nu_n\}$ of measures in $M(X)$
is Cauchy in the metric H.
It follows from the definition of $H$
that $\nu_n(f)$ is Cauchy for every function $f\in\lip(X)$.
By Lemma 4.6, there is a separable measure $\mu$
such that $\nu_n\Rightarrow\mu$.
Show that $H(\nu_n,\mu)\to 0$. Fix $\varepsilon>0$.
Since $\{\nu_n\}$ is Cauchy, there is a number n such
that $H(\nu_n,\nu_m)<\varepsilon$ for all $m>n$.
If $f\in\lip_1(X)$ is some nonnegative Lipschitz function
then the sequence of the bounded functions $f_k(x)=\min\{f(x),k\}$
is monotonically convergent: $f_k\uparrow f$ on $X$.
By Lebesgue's theorem \cite[Chapter 1, \S 12, (12.6), p. 48]{Sax},
$\nu_n(f_k)\uparrow\nu_n(f)$ and $\mu(f_k)\uparrow\mu(f)$;
consequently, there is a number $k$ so large that
$|\nu_n(f)-\nu_n(f_k)|<\varepsilon$ and $|\mu(f)-\mu(f_k)|<\varepsilon$.
It follows from weak convergence that $\lim_{m\to\infty}\nu_m(f_k)=\mu(f_k)$.
By the choice of $n$, we have
$|\nu_n(f_k)-\mu(f_k)|=\lim_{m\to\infty}|\nu_n(f_k)-\nu_m(f_k)|\le\varepsilon$.
Eventually, we obtain
$$|\nu_n(f)-\mu(f)|\le|\nu_n(f)-\nu_n(f_k)|+|\nu_n(f_k)-\mu(f_k)|+|\mu(f_k)-\mu(f)|<3\varepsilon$$
Every Lipschitz function is representable as the difference of two nonnegative functions: $f=f^{+}-f^{-}$;
consequently, $\mu(f)=\lim_{n\to\infty}\nu_n(f)$ for all ${f\in\lip_1(X)}$.
Hence, first, $\mu\in M(X)$, since
$\mu(\phi_{x_0})=\lim_{n\to\infty}\nu_n(\phi_{x_0})<\infty$,
and second, $\lim_{n\to\infty}H(\nu_n,\mu)=0$.
This proves completeness of $M(X)$.
The first part of the theorem is proven.

Assume now that $M(X)$ is complete and show that so is $X$.
Consider some Cauchy sequence
$\{x_n\}_{n=1}^{\infty}$ in $X$.
Since $H(\delta_x, \delta_y)=\rho(x,y)$ for all ${x,y\in X}$,
the sequence $\{\delta_{x_n}\}$ of Dirac measures is Cauchy and,
by completeness of $M(X)$, converges in the metric H to some measure $\nu$.
From separability of $\nu$ we find that the support of $\nu$ is nonempty.
Let $x_0\in\spt\nu$.
Choose an arbitrary $\varepsilon>0$. For the Lipschitz function
$f_{\varepsilon}(x)=\max\bigl(0, 1-\frac{1}{\varepsilon}\rho(x_0,x) \bigr)$
we have $\lim_{n\to\infty}\delta_{x_n}(f)=\nu(f)>0$.
Therefore, there is $N$ such that
$\delta_{x_n}(f)=f(x_n)>0$ for $n>N$; hence,
$\rho(x_0,x_n)<\varepsilon$. Consequently, $x_n\to x_0$ as $n\to\infty$.
Completeness of $X$ is proven.
\hbox{$\Box$}

\vspace{2em}
{\bf 5. Applications of the completeness theorem.} 
As mentioned in the introduction,
Theorem 4.2 has applications in the theory of self-similar sets.
For example, the proof of Hutchinson's theorem \cite[4.4(1)]{Ht}
on existence of a unique invariant measure on an invariant set is
incorrect without Theorem 4.2.
Since regular Borel outer measures are used in \cite[4.4(1)]{Ht},
we should reformulate Theorem 4.2 in terms of outer measures.

An \emph{outer measure} on $X$ is a function
$\mu:2^X\to[0,\infty]$, such that $\mu(\varnothing)=0$ and
$\mu \bigl( \cup_{i=1}^{\infty}E_i \bigr) \le \sum_{i=1}^{\infty}\mu(E_i)$ for $E_i\subset X$.
A set A is called $\mu$-{\it measurable} if $\mu(T)=\mu(T\cap A)+\mu(T\setminus A)$ for all $T\subset X$.
We say that an outer measure $\mu$ is {\it regular Borel} if all
Borel sets are $\mu$-measurable and, for every $A\subset X$, there is
a Borel set $B\supset A$ such that $\mu(A)=\mu(B)$.
We denote by $M^*(X)$ the space of all separable outer measures on $X$
satisfying $\mu(X)=1$ and $\int_X\phi_{x_0}\, d\mu<\infty$ for some $x_0\in X$.
The metric H on $M^*(X)$ is also defined by (1).

\vspace{1em} {\bf Theorem 5.1.} \sl
The space $M^*(X)$ is complete in the metric H if and only if so is $X$. \rm
\vspace{1em} 

The proof is immediate from Theorem 4.2 and the following

\vspace{1em} {\bf Proposition 5.2.} \sl
 The space $M^*(X)$ with the metric H is isometric to $M(X)$. \rm

\vspace{1em} {\tt Proof.}
Let $F:M^*(X)\to M(X)$ act by the formula $F(\mu)=\mu|_{{\cal B}(X)}$.
The measure $F(\mu)$ coincides
with $\mu$ on Borel sets.
It follows from the definition of the Lebesgue integral that the integrals
$\int_X f\,d\mu$ and $\int_X f\,dF(\mu)$
coincide at continuous functions.
Consequently, $F$ preserves the metric $H$.
Since a regular Borel outer measure is determined uniquely
by its values at Borel sets, the mapping $F$ is injective.
For an arbitrary measure $\nu\in M(X)$ the formula
$\mu(A)=\inf\{\nu(B): B\supset A,\, B\in{\cal B}(X)\}$
defines an outer measure on $X$.
The measure $\mu$ is regular Borel,
since by construction for every $A\subset X$,
there is a sequence $B_n\in{\cal B}(X)$, $n=1,2,\ldots$,
such that $B_n\supset A$ and $|\mu(A)-\mu(B_n)|< 1/2^n$,
whence $B=\cap B_n\in{\cal B}(X)$, $B\supset A$, and $\mu(A)=\mu(B)$.
The measure $\mu$ coincides with $\nu$ on Borel sets; consequently,
$\mu\in M^*(X)$ and $F(\mu)=\nu$.
Thus, $F$ is also surjective and hence a one-to-one mapping.
\hbox{$\Box$}

\vspace{1ex}
In \cite{Ak} there is an example of construction of measures in
$\mathbb{R}^n$ invariant under countable systems of
contractions known as IIFS (Infinite Iterated Function Systems) [10]. The following theorem generalizes
this example to the case of a complete metric space.

\filbreak 
\vspace{1em} {\bf Theorem 5.3.} \sl 
In a complete metric space $(X,\rho)$, for every countable system
${\bf S}=\{S_i\}_{i\in \mathbb{N}}$
of contractions
($S_i:X\to X$ and ${\Lip S_i<1}$ for $i\in \mathbb{N}$)
with fixed points xi and for every probability vector
${{\bf p}=(p_1, p_2, \ldots)}$
($\sum_{i\in \mathbb{N}} p_i=1$ and $p_i\ge 0$ for $i\in \mathbb{N}$),
satisfying the condition
$\sum_{i\in \mathbb{N}} p_i \rho(x_1, x_i)<\infty$,
there is a unique measure $\nu\in M(X)$ such that
$$\nu(A)=\sum_{i=1}^{\infty}p_i\nu S_i^{-1}(A)~\mbox{for all}~ \nu-\mbox{measurable}~ A\subseteq X$$
\rm 

{\tt Proof.}
The operator $T:M(X)\to M(X)$, $T(\nu)=\sum_{i=1}^{\infty}p_i\nu S_i^{-1}$,
is a contraction in the metric $H$ (see \cite[3, Theorem 5]{Ak}).
By Theorem 4.2, the space $M(X)$ is complete. Applying the Banach Fixed Point
Theorem, we find that there is a unique measure $\nu\in M(X)$ such that $T(\nu)=\nu$.
\hbox{$\Box$}

\vspace{1em}
The measure $\nu$ generated by the system ${\bf S}$ is called an invariant measure. We can take spt $\nu$ as a set
invariant under ${\bf S}$. Note that the so-constructed sets can be unbounded and thereby noncompact unlike
attractors of iterated function systems. The following example demonstrates also existence of invariant
measures with bounded but not compact supports.
\vspace{1em} \\
{\bf Example 5.1.}
In the Hilbert space $l_2$ with the Hilbert basis $\{e_i\}_{i=1}^{\infty}$,
consider the system ${\bf S}=\{S_i\}_{i=1}^{\infty}$ of contractions
$S_i(x)=\frac{1}{2}(x-e_i)+e_i$
and the probability vector
${\bf p}=\left(\frac{1}{2^i}\right)_{i=1}^{\infty}$.
The support of the corresponding invariant measure is bounded
(lies in the unit ball) but is not compact, since it contains
$\{e_i\}$ that are fixed points of the mappings $\{S_i\}$.

\vspace{1em}
The results of the present article were announced in [11, 12].

\vspace{1ex} 
The author is deeply grateful to his research advisor Professor V.~V.~Aseev for constant attention
to the work and useful pieces of advice during the preparation of the manuscript. The author is also
grateful to the referee for his/her remarks that improved the exposition.

\vspace{1ex}
Author thanks Radu Miculescu from Transilvania University of Brasov, Romania,
for pointing out the flaw in the earlier version of the proof of the Corollary 4.4 \cite[p. 72]{Kr3}.



\indent\hspace{6pt}\sc A. S. Kravchenko, Russia, 2024 \\
\indent\hspace{6pt}\it E-mail address: \tt aleksey14@gmail.com
\end{document}